\documentclass{gen-p-l}

\usepackage{epsf,amssymb}

\newtheorem{thm}{Theorem}[section]
\newtheorem{lem}[thm]{Lemma}
\newtheorem{cor}[thm]{Corollary}

\newtheorem{rem}[thm]{Remark}

\newtheorem*{claim}{Claim}
\newtheoremstyle{definition}{7pt plus6.3pt minus6.3pt}{7pt plus3pt minus3pt}%
{\rm}{}{\bf}{}{0.75em}{\thmname{#1}\thmnumber{ #2}\thmnote{\sl\stdspace#3}}
\theoremstyle{definition}\newtheorem{example}[thm]{Example}
\newtheorem{exercise}[thm]{\small Exercise}

\newcommand{\bbr}{\begin{rem}\em} 
\newcommand{\eer}{\end{rem}}

\newcommand{\bex}{\begin{example}} 
\newcommand{\eex}{\end{example}}

\newcommand{\bhw}{\begin{exercise}\small} 
\newcommand{\ehw}{\end{exercise}}

\newcommand{\be}{\begin{enumerate}}
\newcommand{\ee}{\end{enumerate}}

\def\tb{\operatorname{tb}}
\def\C{\hbox{$\mathbb C$} }
\def\Z{\hbox{$\mathbb Z$} }

\def\R{\hbox{$\mathbb R$} }

\def\co{\colon\thinspace}

\def\dfn#1{{\em #1}}

\begin{document}

\title{Introductory Lectures on Contact Geometry}
\date{June 15, 2000}

\author{John B. Etnyre}
\address{University of Pennsylvania, Philadelphia, PA 19104}
\email{etnyre@math.upenn.edu}
\urladdr{http://www.math.upenn.edu/\char126 etnyre}

\thanks{Supported in part by NSF Grant \# DMS-9705949.}


\maketitle


\section{Introduction}

Though contact topology was born over two centuries ago, in the work of Huygens, Hamilton and 
Jacobi on geometric optics, and been studied by many great mathematicians, such as Sophus Lie, 
Elie Cartan and Darboux, it has only recently moved into the foreground of mathematics.
The last decade has witnessed many remarkable breakthroughs in contact topology, resulting
in a beautiful theory with many potential applications. More specifically, as a coherent -- though
sketchy -- picture of contact topology has been developed, a surprisingly
subtle relationship arose between contact structures and 3- (and 4-) dimensional topology. 
In addition, the applications of contact topology have extended far beyond geometric optics
to include non-holonomic dynamics, thermodynamics and more recently Hamiltonian dynamics
\cite{Hofer93, Weinstein79} and
hydrodynamics \cite{EtnyreGhrist}.

Despite it long history and all the recent work in contact geometry, it is not overly accessible to
those trying to get into the field for the first time. There are a few books giving a brief introduction
to the more geometric aspects of the theory. Most notably the last chapter in \cite{a},
part of Chapter~3 in \cite{McDuffSalamon} and an appendix to the book \cite{Arnold}.
There have not, however, been many books or survey articles
(with the notable exception of \cite{Giroux93}) giving an introduction to the
more topological aspects of contact geometry. It is this topological approach that has lead to
many of the recent breakthroughs in contact geometry and to which this paper is devoted. 
I planned these lectures when asked to give an introduction to contact geometry at the
Georgia International Topology Conference in the summer of 2001. My idea was to give an introduction to
the ``classical'' theory of contact topology, in which the {\em characteristic foliation} plays a
central roll, followed by a hint at the more modern trends, where specific foliations take a back
seat to {\em dividing curves}. This was much too ambitious for the approximately one and a half
hours I had for these lectures, but I nonetheless decided to follow this outline in preparing these
lecture notes. These notes begin with an introduction to contact structures in Section~\ref{dae},
here all the basic definitions are given and many examples are discussed. In the following section
we consider contact structures near a point and near a surface. It is in this section that the 
fundamental notion of characteristic foliation on a surface first appears. 
In an appendix to Section~\ref{ls},
I briefly describe Moser's method, which is a technique for understanding families of contact structures.
Section~\ref{totcs} is devoted to the all pervasive dichotomy in contact geometry: tight vs.\ overtwisted.
Here we see that overtwisted contact structures are not so interesting from a topological point of view and
that tight contact structures have and intimate and subtle relationship with topology. Then,
in  Section~\ref{knots}, we consider
special knots in contact structures. The study of these knots sheds light on 
the tight vs.\ overtwisted dichotomy and allows us to prove a general existence theorem for contact
structures. We end with a brief introduction to convex surfaces. Though this section is short we
will be able to indicate the power of convex surfaces in contact geometry and 
point the interested reader to recent literature on the subject.

These lectures are written in an informal style with many exercises, which are usually
not too difficult and copious hints are provided. 
The proofs of most results are left to the exercises, however for more complicated proofs the outline is
given with the details left as exercises.
I am assuming the reader is familiar with
basic differential topology (manifolds, vector fields, Lie derivatives, forms, $\ldots$, see 
\cite{Warner}) and has a passing knowledge of 3--manifold topology (as can be gleaned from a 
glance or two at \cite{Rolfson} or \cite{Hempel}).

As these notes have a bias for topological techniques in contact geometry, many exciting and
important recent developments have been left out, specifically in regards to the use of
holomorphic curves in contact geometry. Here we refer the reader to \cite{Eliashberg90a, EGH, Gromov85}.
For connections with Seiberg-Witten Theory see \cite{KronheimerMrowka97, Lisca98}. 
Finally for an interesting
historical overview the reader should consult \cite{Geiges01}.

\section{Definitions and Examples}\label{dae}

A \dfn{plane field} $\xi$ on
$M$ is a subbundle of the tangent bundle $TM$ such that $\xi_p=T_pM\cap \xi$ is a 2-dimensional
subspace of $T_pM$ for each $p\in M.$ 
\bex
Consider the 3-manifold $M=\Sigma\times S^1$ where $\Sigma$ is a surface. Then for each $p=(x,\theta)\in
\Sigma\times S^1$ let $\xi_p= T_x\Sigma\subset T_pM.$ Clearly $\xi$ is a plane field on $M.$
\eex
\bex
Let $\alpha$ be a 1-form on $M.$ So at each point $p\in M$  we have a linear map
\begin{equation}
	\alpha_p\co T_pM\to \R.
\end{equation}
Thus $\ker \alpha_p$ is either a plane or all of $T_pM.$ If we assume the 1-form never has
all of $T_pM$ as its kernel, then $\xi=\ker \alpha$ is a plane field. Note in the previous
example the 1-form $\alpha=d\theta$ defines $\xi.$
\eex
It turns out that, locally,  you can always represent a plane field as the kernel of a 1-form.
\bhw
Prove this. In other words, given a plane field $\xi$ on $M$ and a point $p\in M$ show you can
find a neighborhood $U$ of $p$ and a 1-form $\alpha_U$ defined on the neighborhood such that
$\xi|_U=\ker \alpha_U.$
\ehw
\bhw
If $M$ and $\xi$ are both oriented show that you can find a 1-form $\alpha$ defined on all of $M$ such
that $\xi=\ker\alpha.$
\ehw

A plane field $\xi$ is called a \dfn{contact structure} if for any 1-form $\alpha$ with $\xi=\ker\alpha$
($\alpha$ can be locally or globally defined) we have
\begin{equation}
	\alpha\wedge d\alpha\not=0.
\end{equation}
\bhw
Show that $\alpha\wedge d\alpha\not=0$ if and only if $d\alpha\vert_\xi \not=0.$
\ehw
Before we look at some examples of contact structures note that our first example of
a plane field is not a contact structure. Indeed the plane field is defined by the 1-form
$\alpha=d\theta$ so $d\alpha=d(d\theta)=0.$ 

\begin{figure}[ht]
	{\epsfysize=2truein\centerline{\epsfbox{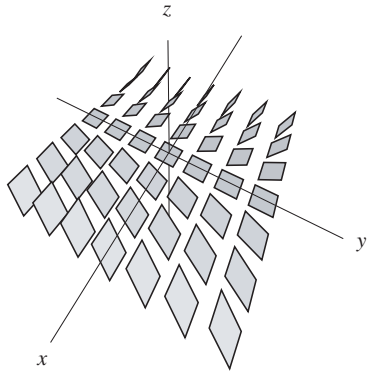}}}
	\caption{The contact structure $\ker(dz+xdy).$}
	\label{ex1}
\end{figure}

\bex\label{standardxi}
Consider the manifold $\R^3$ with standard Cartesian coordinates $(x,y,z)$ and the 1-form
\begin{equation}
	\alpha_1=dz+xdy.
\end{equation}
Note that $d\alpha_1=dx\wedge dy$ so $\alpha_1\wedge d\alpha_1=dz\wedge dx\wedge dy\not=0.$ Thus
$\alpha_1$ is a contact form and $\xi_1=\ker\alpha_1$ is a contact structure. 
At a point $(x,y,z)$ the contact plane $\xi_1$ is spanned by $\{ \frac{\partial}{\partial x},
x\frac{\partial}{\partial z}-\frac{\partial}{\partial y}\}.$ So at any point in the $yz$-plane 
({\em i.e.} where $x=0$) $\xi_1$ is horizontal. If we move to the point $(1,0,0)$ then $\xi_1$ is
spanned by $\{ \frac{\partial}{\partial x},\frac{\partial}{\partial z}-\frac{\partial}{\partial y}\}.$
So the plane is tangent to the $x$-axis but has been tilted clockwise by $45\%.$ In general, if
we start at $(0,0,0)$ we have a horizontal plane and as we move out along the $x$-axis the plane
will twist in a left handed manner ({\em i.e.} clockwise). The twist will be by $90\%$ when $x$ 
``gets to'' $\infty.$ There is similar behavior on all rays perpendicular to the $yz$-plane.
See Figure~\ref{ex1}.

\bbr\label{xysame}
Many authors prefer to use the form $dz-ydx$ to define the ``standard'' contact structure on
$\R^3.$ There is really no difference between these structures. (Rotating about the $z$-axis will
take one of these structures to the other.)
\eer

\eex
\bex\label{symxi}
Consider $\R^3$ with cylindrical coordinates $(r,\theta,z)$ and the 1-form 
\begin{equation}
	\alpha_2=dz+r^2d\theta.
\end{equation}
Since $\alpha_2\wedge d\alpha_2=2rdr\wedge d\theta\wedge dz\not=0,$ $\xi_2=\ker\alpha_2$ is a 
contact structure.
At the point $(r,\theta,z)$ the contact plane $\xi_2$ is spanned by $\{  \frac{\partial}{\partial r},
r^2\frac{\partial}{\partial z}-\frac{\partial}{\partial \theta}\}.$ So when $r=0$ ({\em i.e.} 
in the $z$-axis)
$\xi_2$ is horizontal. As you move out on any ray perpendicular to the $z$-axis the planes $\xi_2$ will twist
in a clockwise manner. So this example is just like the previous one except that everything is 
symmetric about the
$z$-axis.
\eex 

Two contact structures $\xi_0$ and $\xi_1$ on a manifold $M$ are called \dfn{contactomorphic} if there is a 
diffeomorphism $f\co M\to M$ such that $f$ send $\xi_0$ to $\xi_1:$ 
$$f_*(\xi_0)=\xi_1.$$
\bhw
Show that a diffeomorphism $f\co M\to M$ is a contactomorphism if and only if there are contact 
forms $\alpha_0$ and $\alpha_1$ 
for $\xi_0$ and $\xi_1,$ respectively, and a non-zero function $g\co M\to \R$ such 
that $f^*\alpha_1=g\alpha_0.$
\ehw

\bhw
Check that Examples~\ref{standardxi} and \ref{symxi} are contactomorphic. If you are having trouble
coming up with the contactomorphism then first try to write down the contactomorphism implied in 
Remark~\ref{xysame} 
\ehw 

\bex\label{otexample}
Once again consider $\R^3$ with cylindrical coordinates, but this time take the 1-form 
$\alpha_3=\cos r dz +r\sin r d\theta.$ One may compute that 
\begin{equation}
\alpha_3 \wedge d\alpha_3=(1+\frac{\sin r\cos r}{r})d\hbox{vol}.
\end{equation} 
Thus to see that $\alpha_3$ is a contact form you only
have to check that 
\begin{equation}
1+\frac{\sin r\cos r}{r}>0.
\end{equation}
Note that $\xi_3=\ker\alpha_3$ is
horizontal along the $z$-axis and as you move out on any ray perpendicular to the $z$-axis
the planes will twist in a clockwise manner. This time, however, the planes will twist $90\%$ by
the time you get to $r=\pi/2.$ In fact, as you move out on any ray $\xi_3$ will make infinitely
many full twists as $r$ goes to $\infty!$
\eex
This example certainly looks different from our previous two examples, but it is not exactly obvious
how one would actually show it is different. In the early 1980's Bennequin \cite{Bennequin83} 
did distinguish
this example from the previous ones and in the process ushered in a new era in contact geometry. We
will indicate Bennequin's proof in Section~\ref{knots}.

So far all our examples are on $\R^3.$ We now give an example on a closed manifold.
\bex\label{sphere}
Consider the unit 3-sphere, $S^3,$ in $\R^4.$ Let 
\begin{equation}\alpha=(x_1dy_1-y_1dx_1+x_2dy_2-y_2dx_2)|_{S^3},\end{equation}
where $(x_1,y_1,x_2,y_2)$ are standard Cartesian coordinates on $\R^4$ and set $\xi=\ker\alpha.$
\bhw
Check that $\alpha\wedge d\alpha\not=0$ and thus $\xi$ is a contact structure on $S^3.$ \hfill\break
Hint: It
might be helpful to read the following paragraph before trying attempting this exercise. 
\ehw
In anticipation of the next example it will be useful to  describe $\xi$ in another way. If we let
$f(x_1,y_1,x_2,y_2)=x_1^2+y_1^2+x_2^2+y_2^2$ then $S^3=f^{-1}(1).$ Moreover at a point
$(x_1,y_1,x_2,y_2)$ in $S^3$ the tangent space is given by
\begin{equation}
T_{(x_1,y_1,x_2,y_2)}S^3=\ker df_{(x_1,y_1,x_2,y_2)}= \ker (2x_1dx_1+2y_1dy_1+2x_2dx_2+2y_2dy_2).
\end{equation}
Now we can think of $\R^4$ as $\C^2.$ Under this identification we denote the complex 
structure ({\em i.e.} multiplication
by $i$) by $J.$ In other words, $Jx_i=y_i,$
$Jy_i=-x_i$ for $i=1,2.$ The complex structure $J$ induces a complex structure on each tangent 
space: $J\frac{\partial}{\partial x_i}=\frac{\partial}{\partial y_i}$ and 
$J\frac{\partial}{\partial y_i}=-\frac{\partial}{\partial x_i}$ for $i=1,2.$ 
\begin{claim}
The plane field $\xi$ is the set of complex tangencies to $S^3.$ By this we mean 
\begin{equation}
	\xi=T_{(x_1,y_1,x_2,y_2)}S^3\cap J(T_{(x_1,y_1,x_2,y_2)}S^3).
\end{equation}
\end{claim}
Indeed one may easily check that 
\begin{equation}J(T_{(x_1,y_1,x_2,y_2)}S^3)= \ker (df_{(x_1,y_1,x_2,y_2)}\circ J)
\end{equation}
and
\begin{equation}
	df_{(x_1,y_1,x_2,y_2)}\circ J = 2x_1dy_1-2 y_1dx_1+2x_2dy_2-2y_2dx_2.
\end{equation}
Thus we have $\alpha=(df\circ J)|_{S^3}$ and the claim is proved.
\eex

\bhw
Show that $(S^3\setminus\{p\},\xi|_{S^3\setminus\{p\}})$ is contactomorphic to $(\R^3,\xi_2).$
\hfill\break
Hint: Pick the point $p$ carefully and use stereographic coordinates.
\ehw

It turns out that many contact structures can be described as the set of complex tangencies to
a real hypersurface in a complex manifold.
\bex\label{stein}
Let $X$ be a complex manifold with boundary and denote the induced complex structure on $TX$ by $J.$ 
We can find a function $\phi$ defined in a neighborhood of the boundary such that $\phi^{-1}(0)=\partial X.$
Now as in the previous example we can see that the complex tangencies to $M=\partial X$ are 
given by $\ker(d\phi\circ J).$ Thus the complex tangencies $\xi$ to $M$ form a contact structure if
and only if $d (d\phi\circ J)$ is a non-degenerate 2-form on $\xi.$

A fruitful way to construct such manifolds has been through the use of Stein surfaces. To define
Stein surfaces we need some preliminary notions. 
Let $X$ be a complex manifold of complex dimension 2 (real
dimension 4). Again let $J$ denote the induced complex structure on $TX.$ 
From a function
$\phi\co X\to \R$ we can define a 2-form $\omega= d (d \phi \circ J)$ and a symmetric form
$g(v,w)=\omega(v, Jw).$ If this symmetric form is positive definite ({\em i.e.} defines a metric
on $X$) the function $\phi$ is called \dfn{strictly plurisubharmonic}. The manifold 
$X$ is a \dfn{Stein surface} if $X$ admits
a proper strictly plurisubharmonic function $\phi\co X\to \R.$ It is easy to see that in this situation
the complex tangencies to $M_c=\phi^{-1}(c)$ form a contact structure whenever $c$ is not a
critical value. 
We will call such a contact structure
\dfn{Stein fillable}. Later we will see that this implies $\xi$ is a special type of contact structure.
See \cite{Gompf98} to learn how to construct many Stein surfaces and hence many contact structures.
\eex

\section{Local Structure}\label{localstructure}\label{ls}

In this section 
we discuss the nature of contact structures near a point (Darboux's Theorem) and near a 
surface. You can find further discussion of all these local theorems in \cite{a, McDuffSalamon}.

\subsection{Darboux's Theorem}

Darboux's Theorem essentially says that all contact structures look the same near a point. So
contact structures do not have interesting local structure (this should be compared with
Riemannian geometry, where the curvature is an obstruction to metrics being locally the same).
This is an indication that any interesting phenomena in contact geometry should be of a global
nature ({\em i.e.} be related to the global topology of the manifold supporting the contact structure).

\begin{thm}\label{darboux}
	Let $(M,\xi)$ be any contact 3-manifold and $p$ any point in $M.$ Then there
	exist neighborhoods $N$ of $p$ in $M,$ and $U$ of $(0,0,0)$ in $\R^3$ and
	a contactomorphism 
	$$f\co(N,\xi|_N)\to(U,\xi_1|_U).$$
\end{thm}

The current modern proof of Darboux's Theorem uses ``Moser's Method.'' We will discuss this more in 
the appendix to this section. The classical proofs of this theorem are more elementary in nature and
we encourage the reader to try and come up with an elementary proof.
\bhw
Find an elementary proof of Darboux's theorem.
\ehw

\subsection{The Characteristic Foliations}

\begin{figure}[ht]
	{\epsfysize=1.2truein\centerline{\epsfbox[0 719 130 841]{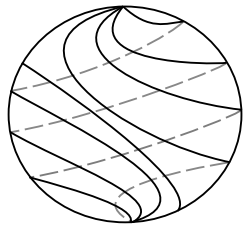}}}
	\caption{The characteristic foliation on  $S^2.$}
	\label{charfolonS3}
\end{figure}

Let $\Sigma$ be an embedded oriented surface in a contact manifold $(M,\xi).$ 
At each point $x$ of $\Sigma$ consider 
$$l_x=\xi_x\cap T_x\Sigma.$$
For most $x,$ the subspace $l_x$ will be a line in $T_x\Sigma,$ but at some points, which we call
singular points, $l_x=T_x\Sigma.$
\bhw
Show that $l_x$ cannot equal $T_x\Sigma$ for all $x$ in some open subset of $\Sigma.$ \hfill\break
Hint: If this is true, then you can show that the contact condition is violated.
Consider two vectors fields $v$ and $w$ tangent to $\Sigma$ defined along this open subset.
Using the formula $d\alpha(v,w)=v\alpha(w)-w\alpha(v)+\alpha([v,w])$ compute $\alpha\wedge d\alpha.$
\ehw
It is not hard to show (see the next exercise) that
we may find a singular foliation $\mathcal{F}$ of $\Sigma$ tangent to $l_x$ at each $x.$ By this
we mean the complement of the singularities is the disjoint union of 1-manifolds, called leaves 
of $\mathcal{F},$ and the leaf through $x$ is tangent to $l_x.$ This singular foliation is called
the \dfn{characteristic foliation} of $\Sigma$ and is 
denoted $\Sigma_\xi$ (some authors prefer $\xi\Sigma$).
\bhw
Show there is a singular foliation tangent to $l_x.$\hfill\break
Hint: Locally on the surface find a vector field tangent to $l_x$ at the nonsingular points and
$0$ at the singular points. Then use basic existence results from ordinary differential equations
to construct the singular foliation locally. Finally make sure your local foliations fit together
to give a global foliation. 
\ehw
\bex
Let $\Sigma$ be the unit sphere in $(\R^3,\xi_2).$ The only singularities in $\Sigma_{\xi_2}$ are at
the north and south poles. See Figure~\ref{charfolonS3}.
\eex

\begin{figure}[ht]
	{\epsfysize=1.5truein\centerline{\epsfbox[148 337 455 464]{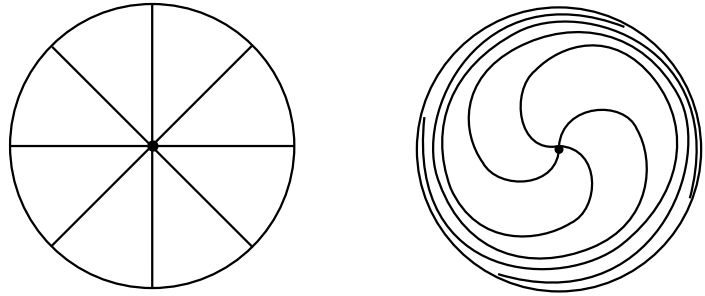}}}
	\caption{The degenerate and non-degenerate characteristic foliation on $D^2.$}
	\label{diskfol}
\end{figure}

\bex
Let $\Sigma$ be the disk of radius $\pi$ in the $r\theta$-plane in $(\R^3,\xi_3).$ 
As shown on the left hand side of in Figure~\ref{diskfol}, the center of
$\Sigma$ is a singular point and each point on the boundary of $\Sigma$ is also a singular point.
Let $\Sigma'$ be $\Sigma$ with its interior pushed up slightly. 
Now the only singularity in the characteristic foliation is
at the center point. The boundary of $\Sigma'$ is now a closed leaf in the foliation. See 
Figure~\ref{diskfol}.
\eex

\begin{figure}[ht]
	{\epsfysize=1.5truein\centerline{\epsfbox[152 345 472 470]{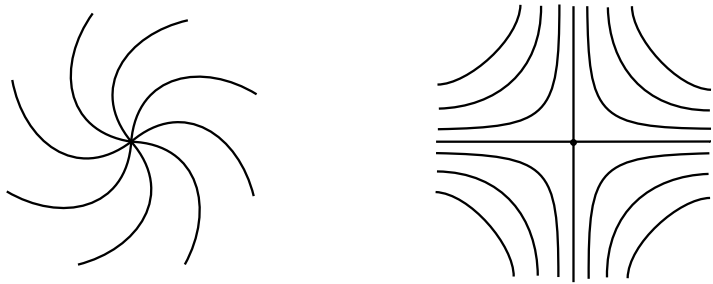}}}
	\caption{Generic singularities in the characteristic foliation.}
	\label{sings}
\end{figure}

This last example illustrates an important point: any surface $\Sigma$ may be perturbed by a 
$C^\infty$-small isotopy so that its characteristic foliation has only ``generic'' isolated singularities.
A singularity is ``generic'' if it looks like one in Figure~\ref{sings}. On the left hand side is an
elliptic point and on the right hand side is a hyperbolic point.

Recall that $\Sigma$ is oriented and we chose an orientation on $\xi.$ We can orient 
$l_x=\xi_x\cap T_x\Sigma$ as follows: the vector $v\in l_x$ orients $l_x$ if when we choose  vectors
$v_\xi\in \xi_x$ and $v_\Sigma\in T_x\Sigma$ such that $(v,v_\xi)$ orients $\xi_x$ and $(v,v_\Sigma)$ orients
$T_x\Sigma$ then $(v,v_\xi,v_\Sigma)$ orients $M.$   In this manner the leaves of the characteristic
foliation inherit an orientation, so we can draw arrows on the leaves of the foliation and think of
the foliation as a ``flow.'' Moreover, at each singular point $x$ we can assign a sign: the singularity
is $+$ (respectively, $-$) if the orientation on $T_x\Sigma$ agrees 
(respectively, disagrees) with the one on $\xi_x.$ With these
conventions, a positive elliptic point is a source, a negative elliptic point is a sink.
Note the sign if a hyperbolic point is not obvious at first glance as it is related to the
ratio of the eigenvalues associated to the linearized flow at the singularity.
\bhw
Determine the sign of a hyperbolic point.
\ehw

\begin{thm}\label{foldet}
	Let $(M_i,\xi_i)$ be a contact manifold and $\Sigma_i$ an embedded surface for $i=0,1.$
	If there is a diffeomorphism $f\co\Sigma_0\to\Sigma_1$ that preserves the characteristic
	foliation:
	$$f((\Sigma_0)_{\xi_0})=(\Sigma_1)_{\xi_1},$$
	then $f$ may be extended to a contactomorphism in some neighborhood of $\Sigma_0.$
	Moreover, if $f$ was already defined on a neighborhood of $\Sigma_0$ then we can
	isotop $f$ so as to be a contactomorphism in some (possibly) smaller neighborhood.
\end{thm}
So the characteristic foliation of on a surface determines (the germ of) the contact structure
near the surface. Once again this theorem my be proved using ``Moser's method''.

\medskip
\begin{center}{\sc Appendix to Section~\ref{localstructure}: Moser's Method}
\end{center}

There are several good references for Moser's method and its many corollaries \cite{a, McDuffSalamon}.
One of the most general theorems one can prove using these techniques is
\begin{thm}
	Let $M$ be an oriented three manifold and $N\subset M$ a compact subset. Suppose
	$\xi_0$ and $\xi_1$ are contact structures on $M$ for which $\xi_0|_N=\xi_1|_N.$ Then 
	there is a neighborhood $U$ of $N$ such that the identity map on a neighborhood of $N$ 
	is isotopic, rel. N, to a contactomorphism when restricted to $U.$
\end{thm}
\bhw
Show Theorems~\ref{darboux}, \ref{foldet} and \ref{transnbhd} follow from this theorem.\hfill\break
Hint: Consider Darboux's theorem. Write down a diffeomorphism from a neighborhood $N'$ of the point
$p$ in $M$ to a neighborhood $U'$ of $(0,0,0)$ in $\R^3$ so that the contact plane at $p$ is sent
to the contact plane at $(0,0,0).$ Push the contact structure $\xi$ forward to $U.$ Now you have
two contact structures on $U$ that agree on $(0,0,0),$ so use the above theorem to finish the proof. 
Theorem~\ref{transnbhd} is similar to this and
Theorem~\ref{foldet} is similar but it is not so obvious you can write down the correct initial 
diffeomorphism.
\ehw

The proof of this theorem follows essentially from the above mentioned references, but we
will outline the proof in the following exercises.
\bhw
Let $\alpha_i$ be a contact form for $\xi_i, i=0,1,$ that determines the orientation on $\xi_i|_N.$
Let $\alpha_t=(1-t)\alpha_0+t\alpha_1.$ Show that on some neighborhood $U'$ of $N$ all the 
$\xi_t=\ker \alpha_t$'s are 
contact structures.
\ehw
\bhw
We now wish to find a family of diffeomorphisms $\phi_t\co U\to U''$ ($U$ 
and $U''$ are possibly smaller neighborhood
of $N$) such that $\phi_t^*\xi_t=\xi_0.$ (Here $\phi_t^*=(\phi_t^{-1})_*.$) This 
will of course finish the proof of the theorem. We will find the $\phi_t$'s as the flow of a vector field.
Suppose $v_t$ is a time dependent vector field whose flow generates the $\phi_t$'s. Show that 
if $v_t\in \xi_t$ then the
$\phi_t$'s satisfy $ \phi_t^*\xi_t=\xi_0$ if and only if $\iota_{v_t}d\alpha_t|_{\xi_t}=
(h_t\alpha_t-\frac{d\alpha_t}{dt})|_{\xi_t},$ where $h_t=\frac{d\alpha_t}{dt}(X_t)$ and $X_t$ is the unique vector
field satisfying $\alpha_t(X_t)=1$ and $\iota_{X_t}d\alpha_t=0.$ 
(Here $\iota_{v_t}$ means contraction with $v_t$). 
\ehw
\bhw
Given $\alpha_t$ above, prove there is a $v_t$ as described in the previous exercise.
\ehw

\section{Tight and Overtwisted Contact Structures}\label{totcs}

There is a fundamental dichotomy in 3-dimensional contact geometry. A contact structure $\xi$ 
on $M$ is called \dfn{overtwisted} if there is an embedded disk $D$ whose characteristic foliation
is homeomorphic to the either one shown in Figure~\ref{diskfol}.
Such a disk is called an \dfn{overtwisted disk}.
A contact structure is called \dfn{tight} if it does not contain an overtwisted disk. Though tight vs.\
overtwisted is obviously a dichotomy, it is not clear that it is a useful one. Throughout the
rest of these lectures we will indicate that overtwisted contact structures are somewhat ``easy'' to
deal with, whereas tight contact structures are quite a bit more difficult to understand. Moreover,
a tight contact structure is capable of detecting subtle properties of the manifold supporting it.

Later we will see directly that overtwisted contact structures are fairly simple to
construct and work with. This is all reflected in the following theorem.
\begin{thm}[Eliashberg \cite{Eliashberg89}]\label{otclass}
	Given a closed compact 3-manifold $M,$ let
	$\mathcal{H}$ be the set of homotopy classes of (oriented) plane fields on $M$ and 
	$\mathcal{C}_o$ be the set of isotopy classes of (oriented) overtwisted contact structures
	on $M.$ The natural inclusion map $\mathcal{C}_o$ into $\mathcal{H}$ induces a homotopy
	equivalence.
\end{thm}
This theorem basically reduces the classification of overtwisted contact structures on a 3-manifold
to the classification of homotopy classes of plane fields. This latter problem is algebraic in
nature and can be understood through the Thom-Pontryagin construction, see \cite{Milnor}. In addition, see
\cite{Gompf98} for a discussion with contact geometry in mind.

One thing, among many, that this theorem implies is that any 3-manifold has an overtwisted contact
structure on it! Moreover, any $c\in H^2(M,\Z)$ that is the Euler class of an oriented plane field
is also the Euler class of an overtwisted contact structure.
\bhw
Show that $c\in H^2(M,\Z)$ is the Euler class of
an oriented plane field if and only if its mod 2 reduction is $0.$ (You might need to review a few facts
about characteristic classes to do this.)
\ehw
\noindent
Tight contact structures are not understood nearly as well and they do not always exist.
\begin{thm}[Etnyre-Honda \cite{EtnyreHonda01}]
	There exists a closed compact 3-manifold that does not support any
	tight contact structure.
\end{thm}
Despite this theorem, it seems that in some sense ``most'' 3-manifolds do admit tight
contact structures and when they do they reveal interesting things about the manifold, see 
Section~\ref{genusbounds} 
below. The easiest, and most common, way to construct tight contact structures is via
symplectic geometry. Recall a closed two form $\omega$ on a 4--manifold $X$ is a \dfn{symplectic
form} if $\omega\wedge\omega\not=0.$ A compact symplectic 4--manifold $(X,\omega)$ is said
to fill a contact 3--manifold $(M,\xi)$ if $\partial X=M$ (as oriented manifolds!) and
$\omega|_\xi$ is an area form on $\xi.$ Note that all Stein fillable contact structures 
(Example~\ref{stein}) are filled by a symplectic 4--manifold (since $\omega=d(d\phi\circ J)$ is a
symplectic form). 
\begin{thm}[Eliashberg, Gromov \cite{Gromov85, Eliashberg90a}]\label{fill2tight}
	If a contact structure can be filled by a compact symplectic manifold then
	it is tight.
\end{thm}

We will not go into what is known about the classification of tight contact
structures, see \cite{Giroux00, Honda1, Honda2}, but we do mention the method most commonly
used to understand them. The key ingredient in all classification
results is the following:
\begin{thm}[Eliashberg \cite{Eliashberg92}]\label{3ball}
	If $\mathcal{F}$ is a singular foliation on $S^2$ that is induced by some tight
	contact structure, then there is a unique (up to isotopy fixing the boundary) tight
	contact structure $\xi$ on $B^3$ such that $(\partial B^3)_\xi=\mathcal{F}.$
\end{thm}
Now to understand tight contact structures on a manifold $M$ one ``merely'' removes pieces from
$M$ on which you understand the contact structure ({\em e.g.} neighborhoods of surfaces on
which the characteristic foliation is known) until all that is left of $M$ is a collection 
of 3-balls. Then apply the previous theorem to conclude you understand the contact structure.
This is, of course, quite vague but to understand the strategy better try the following exercise.
\bhw
By Theorem~\ref{fill2tight} the contact structure on $S^3$ described in Example~\ref{sphere} is tight.
Use the above strategy to show there is only one tight contact structure on $S^3.$ Specifically, fill in
the details and understand the following argument: If you have two tight contact structures on
$S^3$ use Darboux's Theorem to say they agree in a neighborhood of a point. Then use Theorem~\ref{3ball}
to conclude that they agree in the complement of the neighborhood. 
\ehw

\subsection{Manipulations of the Characteristic Foliations}

Since any 3--man\-i\-fold can be cut up along surfaces into a collection of 3--balls (in many
ways, {\em e.g.}\ Heegaard decompositions, Haken decompositions, $\ldots$ ) it is clear, from the
strategy discussed above, that 
to understand tight contact structures on a 3--manifold we should understand tight contact structures
in the neighborhood of surfaces better. A first step in this direction is to develop techniques
to manipulate characteristic foliations. One of the most important theorems along these lines
is:
\begin{lem}[Elimination Lemma: Giroux, Fuchs \cite{Eliashberg93}]\label{elem}
	Suppose $\gamma$ is a leaf in a characteristic foliation $\Sigma_\xi$ connecting
	an elliptic and hyperbolic point of the same sign. Then given any neighborhood $N$ of $\gamma,$ 
	we may find an isotopy, supported
	in $N,$ of $\Sigma$ to $\Sigma'$ so that $\Sigma'_\xi\cap N$ contains no singularities and,
	of course, $\Sigma_\xi$ and $\Sigma'_\xi$ agree outside of $N.$ See Figure~\ref{cancel}.
\begin{figure}[ht]
	{\epsfysize=1.2truein\centerline{\epsfbox[87 682 444 766]{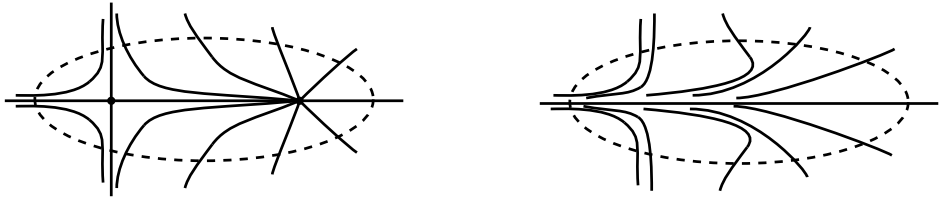}}}
	\caption{The cancellation of singularities with the same sign.}
	\label{cancel}
\end{figure}
\end{lem}
Thus this theorem says we may eliminate singularities of the same sign that are connected by an arc!
\bhw\label{create}
Visualize Figure~\ref{cancel} 
in $(\R^3, \xi_1)$ (recall $\xi_1=\ker{dz+xdy}$) as follows. Start with a embedded rectangle containing
the $y$-axis and tilted slightly out of the $xy$-plane ({\em e.g.} a piece of the graph of 
$f(x,y)=\epsilon x$). The characteristic foliation on this is nonsingular. Now create two singularities
by rotating the middle part of the rectangle past the $xy$-plane ({\em e.g.} rotate a bit of the
rectangle to agree with the graph of $-f(x,y)$). If you did this correctly then the characteristic
foliation should look like the left hand side of Figure~\ref{cancel}. From the construction we know
how to remove the singularities in this example. Use Theorem~\ref{foldet} to prove Lemma~\ref{elem}.
This argument is explicitly worked out in \cite{a}.
\ehw

There is an important strengthening of the Elimination Lemma. Note that in the Elimination Lemma
the arc $\gamma$ is part of some leaf of the new characteristic foliation on $\Sigma'.$ The strengthened
lemma give some control over this new leaf.
\begin{lem}\label{lem:create}
	Suppose $\gamma$ is as in the Elimination Lemma. Let $\gamma'$ be any leaf (distinct from 
	$\gamma$) that limits to the same elliptic point as $\gamma.$ Then we may assume that after 
	the cancellation of the singularities $\gamma$ and $\gamma'$ are on the same leaf of the new
	characteristic foliation.
\end{lem}
Note that there is no flexibility over which two leaves limiting to a hyperbolic point will end up on 
the same leaf after the cancellation.

As Exercise~\ref{create} indicates, it is much easier to create singularities that eliminate them. 
In particular we have
\begin{lem}\label{lem:create2}
	Let $\gamma$ be a segment of a leaf in $\Sigma_\xi$ and $N$ be a neighborhood of $\gamma$
	such that $\Sigma_\xi\cap N$ contains no singularities.
	Then we may find an isotopy, supported in $N,$ of $\Sigma$ to $\Sigma'$ so that $\Sigma'_\xi\cap
	N$ contains an elliptic and hyperbolic singularity of the same sign and $\Sigma_\xi$ and
	$\Sigma'_\xi$ agree outside of $N.$
\end{lem}

Up to this point the careful reader might have been concerned that we discuss ``elliptic''
singularities as if there were only one type of elliptic singularity. (A similar discussion
applies to hyperbolic singularities.) Topologically this
is true ({\em i.e.} up to homeomorphism) but up to diffeomorphism this is not true and Theorem~\ref{foldet}
needs a diffeomorphism! 
\bhw
Show that any two elliptic sources singularities are topologically equivalent (and similarly for sinks).
\hfill\break
Hint: This is a small extension of the Hartman-Grobman Theorem which you can find most books on 
dynamical systems \cite{R}. 
\ehw
\bhw
Show that (generically) up to ($C^1$) diffeomorphism an elliptic singularity is determined by 
the eigenvalues of its linearization (this is not so easy, you might want to consult 
\cite{R}).
\ehw
So how is it that we can ignore this subtlety? It turns out that we may perturb a surface near
an elliptic singularity so that the singularity will be diffeomorphic to a preassigned 
elliptic singularity. 
\bhw
Verify this statement.\hfill\break
Hint: Use Darboux's theorem to reduce the problem to considering disks in $(\R^3, \xi_1)$ which are
tangent to the $xy$-plane at the origin. 
Such disk can be represented as graphs of functions $\{(x,y, f(x,y))\}.$
Now use the previous exercise and perturbations of $f$ to prove the statement. 
\ehw
So as long as we are willing to perturb our surfaces (by a $C^\infty$-small
isotopy) we may ignore this problem of smooth equivalence of elliptic singularities. More precisely,
we actually have
\begin{lem}
	Suppose there is a homeomorphism from $\Sigma_\xi$ and $\Sigma'_{\xi'}$ (both characteristic
	foliations should be generic), then there is a $C^\infty$-small isotopy of $\Sigma'$ to $\Sigma''$
	such that $\Sigma_\xi$ and $\Sigma''_{\xi'}$ are diffeomorphic by a diffeomorphism that
	is isotopic to the original homeomorphism. 
\end{lem}
Thus we can just ``look at'' the characteristic foliation and do not need to worry about the
subtleties of the singularities.

\subsection{Tightness and Genus Bounds}\label{genusbounds}

We use the above manipulations of the characteristic foliation to show 
\begin{thm}[Eliashberg \cite{Eliashberg92}]\label{genusbound}
	Let $(M,\xi)$ be a tight contact 3-manifold and $\Sigma$ an embedded surface in
	$M.$ If $e(\xi)\in H^2(M,\Z)$ denotes the Euler class of $\xi,$ then
	\begin{equation}\label{genusin}
		|e(\xi)([\Sigma])|\leq \begin{cases}
      			\quad -\chi(\Sigma) & \text{ if } \Sigma\not= S^2,\\
        \quad 0 & \text{ if } \Sigma=S^2,
        \end{cases}
	\end{equation}
	where $[\Sigma]$ denotes the homology class of $\Sigma.$
\end{thm}

Though it may not be apparent at first,
this theorem begins to indicate the delicacy of tight contact structures. For example, we have
\begin{cor}
	There are only finitely many elements in $H^2(M,\Z)$ that can be realized as
	the Euler class of a tight contact structure.
\end{cor}
\begin{proof}
	There is no torsion in $H_2(M,\Z).$
	\bhw
		Show this.\hfill\break
		Hint: Use Poincar\'e Duality and the Universal Coefficients Theorem.
	\ehw
	\noindent
	Now let $g_1, \ldots, g_n$ be generators for $H_2(M,\Z)$ 
	\bhw
		Show that any element in $H_2(M,\Z)$ can be represented by a surface.
	\ehw
	\noindent
	Let $S_1, \ldots, S_n$ be embedded surfaces such that the homology class of $S_i$ is $g_i,$
	for $i=1, \dots, n.$ We can assume that none of the $S_i$ are 2-spheres (Why?) and then for each
	of the $S_i,$  Inequality~\eqref{genusin} gives a region between two parallel hyperplanes 
	in $H^2(M,\Z)$ in which an Euler class for a tight contact structure can live.
	\bhw
		Show that all the hyperplanes coming from the $S_i$ define a compact convex
		polytope in $H^2(M,\Z).$
	\ehw 
	\noindent
	There can clearly be only finitely many Euler classes of tight contact structures since they
	have to live in this polytope.
\end{proof}

Note that this corollary clearly shows the difference between tight and overtwisted contact
structures, since any element in $H^2(M,\Z)$ whose mod 2 reduction equals 0 is the Euler class
of an overtwisted contact structure by Theorem~\ref{otclass}.

Inequalities like \eqref{genusin} have shown up in other places too. For example, Thurston 
\cite{Thurston} proved
that the inequality in Theorem~\ref{genusbound} is true for the Euler class of a taut foliation. 
Due in part to this inequality, and many interesting constructions, foliation theory has found
a central place in 3--manifold topology.

\begin{proof}[Proof of Theorem~\ref{genusbound}]
Note that it suffices to prove the theorem when $\Sigma$ is connected.
We begin by trying to understand how to calculate $e(\xi)([\Sigma])$ and $\chi(\Sigma)$ in terms
of $\Sigma_\xi.$ First perturb $\Sigma$ so that the characteristic foliation is generic. (By generic,
we mean that the singularities are isolated elliptic or hyperbolic points and no two hyperbolic
points are connected by a leaf in the foliation.) Then
let $e_\pm$ be the number of $\pm$ elliptic points in $\Sigma_\xi$ and $h_\pm$
be the number of $\pm$ hyperbolic points in $\Sigma_\xi.$  We first have the following simple 
observation:
\begin{equation}\label{eqn:chi}
	\chi(\Sigma)= (e_++e_-)-(h_++h_-).
\end{equation}
This should be clear since we may take a vector field $v$ that directs the characteristic foliation
({\em i.e.} is tangent to $\Sigma_\xi$ at non-singular points, is zero at the singularities and
induces the orientation on $\Sigma_\xi$). The Poincar\'e-Hopf Theorem \cite{Milnor} 
now says that $\chi(\Sigma)$
can be computed in terms of the zeros of $v.$
\bhw
Check that the Poincar\'e-Hopf Theorem implies Equation~\eqref{eqn:chi}.
\ehw 

We now claim that
\begin{equation}\label{eqn:e}
	e(\xi)([\Sigma])= (e_+-h_+)-(e_--h_-).
\end{equation}
To see this recall the the Euler class of a bundle is the obstruction to finding a
non-zero section of the bundle. Moreover, $e(\xi)([\Sigma])$ is just the Euler class
of the restriction of $\xi$ to $\Sigma$ (since everything behaves well with respect to pull
back). Thus to compute the Euler class of $\xi|_\Sigma$ we just need to take a generic 
section of $\xi|_\Sigma$ and calculate the intersection of its graph (in $\xi|_\Sigma$)
with the zero section. 
More specifically, take $v$ from above as our section then the graph of $v$ is 
\[\Gamma=\{(x,p)\in \xi|_\Sigma : p=v(x)\},\]
where $x$ is a point in $\Sigma$ and $p\in \xi_x.$ So $\Gamma$ is a surface in the 4-manifold
$\xi|_\Sigma,$ this is a 4-manifold since it is the total space of a 2-dimensional vector bundle
over a surface. The zero section, $\Gamma_0=\{(x,0)\in \xi|_\Sigma\},$ is another surface. Now
the Euler class of $\xi|_\Sigma$ is just the (oriented) intersection number of these two surfaces. 
\bhw
Show that the contribution to the intersection number of each zero of $v$ is a $+1$ for a positive
elliptic or negative hyperbolic point and a $-1$ for a negative elliptic or positive hyperbolic point.
\hfill\break
Hint: It might be helpful to think about $\chi(\Sigma)$ in these terms and remember $\xi_x=\pm T_x\Sigma$
at the singularities.
\ehw

Now to prove Equation~\eqref{genusin} when $\Sigma\not= S^2$ we need to see that 
$\pm e(\xi)([\Sigma])\leq -\chi(\Sigma).$
Adding Equations~\eqref{eqn:chi} and \eqref{eqn:e} we see that
\begin{equation}
	\chi(\Sigma)+e(\xi)([\Sigma])=2(e_+-h_+).
\end{equation}
So if we can show that, after isotoping $\Sigma,$ $e_+=0,$ then we will know 
$e(\xi)([\Sigma])\leq -\chi(\Sigma).$ To this end, we first arrange that there are no
closed leaves in $\Sigma_\xi$ by using Lemma~\ref{lem:create} to
creating negative elliptic-hyperbolic pairs along any 
closed leaf. (Note we will of course have to isotop $\Sigma$ to do this, but we still call
the resulting surface $\Sigma.$) Now if there are any positive elliptic points $x$ then
let $U_x$ be the set of all leaves in $\Sigma_\xi$ that limit to $x$ and
$B_x$ be the closure of $U_x.$ Denote $B_x\setminus U_x$ by $\partial B_x.$  Ultimately we
will show that $\partial B_x$ contains a positive hyperbolic point $y$ (which is clearly
connected to $x$ by an arc) and use the Elimination
Lemma to cancel $x$ and $y.$  To do this we need to understand the structure of $B_x$ better.

Refer to Figure~\ref{bx} as we discuss $B_x.$ 
First note that $U_x$ does not contain any singularities (other than $x$), 
so all the singularities in $B_x,$ except $x,$ 
are in $\partial B_x.$ 
Now if $p$ is an elliptic point in $\partial B_x$ then it must be negative (Why?). Also note that if $p$ is
a hyperbolic point in $\partial B_x$ then its unstable manifolds (the two curves in $\Sigma_\xi$ that
limit to $p$ in backwards time when we think if $\Sigma_\xi$ as a flow) 
are also in $\partial B_x$ and they limit (in forward time) to negative
elliptic points in $\partial B_x.$
\begin{figure}[ht]
	{\epsfysize=1.6truein\centerline{\epsfbox[98 560 354 770]{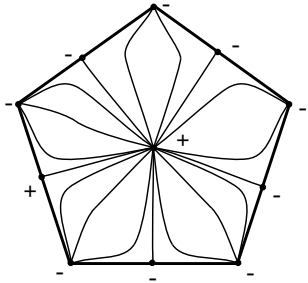}}}
	\caption{A typical $B_x.$}
	\label{bx}
\end{figure}

We claim that $\partial B_x$ contains a positive hyperbolic point. To see this,
we assume that there are no positive hyperbolic points in $\partial B_x$ and 
derive a contradiction. Note that $U_x$ is
embedded in $\Sigma$ and is diffeomorphic to an open disk. 
\bhw
Show that if $B_x$ is embedded then it is diffeomorphic to a closed disk, with piecewise smooth boundary,
and the boundary of the disk
contains only negative elliptic and hyperbolic points connected by arcs.
\ehw
Thus if $B_x$ is embedded then we may use the (strengthened) Elimination Lemma to cancel all the
singularities in $\partial B_x$ resulting in an overtwisted disk. Thus $B_x$ cannot be embedded.

\bhw
If there are hyperbolic points in $\Sigma_\xi,$ then
convince yourself that we can think of $B_x$ as the image of an immersed polygon $f\co P\to M$ such that 
$f$ is an embedding on the interior of $P.$ Moreover, it can be arranged that each edge maps to
the union of a hyperbolic point and its unstable manifolds and each vertex maps to an elliptic point.
See Figure~\ref{bx}. Though $B_x\subset \Sigma$ we will sometimes talk as if $B_x=P,$ this will simplify
notation and (hopefully) clarify what is going on. Just remember that $B_x$ can refer to the image of
an immersed polygon or the polygon itself depending on context.
If there are no hyperbolic singularities in $\Sigma_\xi$ then convince yourself that $\Sigma=S^2,$ $x$ is
the unique positive elliptic point in $\Sigma_\xi$ and $\partial B_x$ is the unique negative elliptic point
in $\Sigma_\xi.$ 
\ehw\label{possibleB}
If $B_x$ is not embedded then $f$ identifies vertices, or vertices and edges, of $P.$ 
Suppose $f$ identifies only vertices.
In this case one may refine 
Lemma~\ref{lem:create2} to create
a negative elliptic--hyperbolic pair near each non-embedded vertex, as shown in Figure~\ref{vert}, 
so as to make $B_x$ embedded for this new characteristic foliation.
\begin{figure}[ht]
	{\epsfysize=1.4truein\centerline{\epsfbox[9 686 400 835]{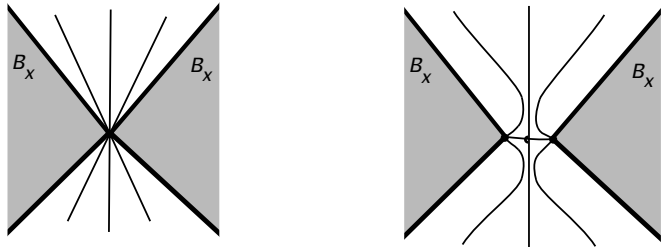}}}
	\caption{Making vertices disjoint.}
	\label{vert}
\end{figure} 
Thus we are back in the embedded case and can construct
an overtwisted disk. 

We are left to consider the case when $f$ identifies edges of $P.$ We consider the simplest case
first. Suppose the image of $f$ is as shown on the left hand side of Figure~\ref{edges}.
 \begin{figure}[ht]
	{\epsfysize=1.4truein\centerline{\epsfbox{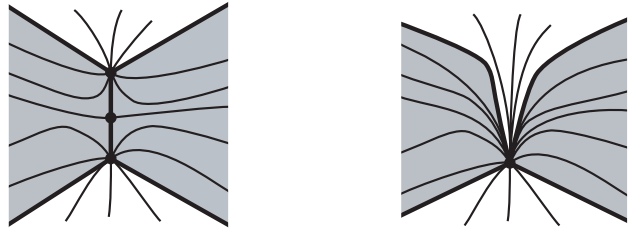}}}
	\caption{A possible $B_x$ when edges are identified (right) and the resulting $B_x$ when
	two singularities are canceled (left).}
	\label{edges}
\end{figure}
If we cancel the hyperbolic point with the upper elliptic point then the new $B_x$ will be
related to the old $B_x$ as shown in Figure~\ref{edges}. Thus the new $B_x$ has only vertices identified, but
we know from here we can get to an embedded $B_x$ and thus an overtwisted disk! 

Note the right hand side
of Figure~\ref{edges} is not correct if the top and bottom vertices on the left hand side are
identified. In this case a periodic obit will be formed. If this happens then we loose the structure
of $B_x.$ To prevent this from happening always cancel edges with distinct vertices first
(making vertices disjoint whenever possible). Since $\partial B_x$ is connected we will
eventually get to the situation where there is only one vertex in $\partial B_x.$ 
If there are no edges left then $\Sigma=S^2$ as we discussed above. If there is only one edge
then $B_x$ is embedded and after canceling the boundary singularities we have an overtwisted 
disk. If $\partial B_x$ has two or more edges that are not identified then we may use the move
depicted in Figure~\ref{vert} to create two distinct vertices on $\partial B_x\subset \Sigma$ and
cancel more edges. So if we have simplified $B_x$ as far as possible and have not found an
overtwisted disk (or $\Sigma=S^2$) then $\partial B_x\subset \Sigma$ has only one vertex and 
all but possibly one edge is identified with some other edge. Thus the image of $B_x$ in $\Sigma$ is a 
closed surface or a subsurface with one boundary component. 
We show how to construct an overtwisted disk when 
$B_x=\Sigma=T^2.$ In this case $\partial B_x\subset \Sigma$ has two edges $E$ and $F.$ 
Let $N$ be a neighborhood of $E$ in $\Sigma$ and $N'$ (respectively $N''$) a copy of 
$N$ pushed slightly up off of $\Sigma$ (respectively slightly down off of $\Sigma$). We can now
cut $\Sigma$ along $E$ and push one side up to agree with $N'$ and the other edge down
to agree with $N''.$  Note since the foliation on $N$ is generic it is also stable so $N'_\xi=
N''_\xi=N_\xi.$ If we consider $B_x$ sitting the new surface then 
we now have two copies of $E$ and two vertices.
\bhw
Find an overtwisted disk associated to $B_x$ on this new surface.
\ehw
Even though this new $B_x$ is not on $\Sigma$ this is not a
problem since we are arguing by contradiction --- using our assumption 
that $\partial B_x$ has no positive hyperbolic singularities to
construct an overtwisted disk in $M.$ It is irrelevant that the overtwisted disk is not actually on
$\Sigma.$   
\bhw
Expanding the above argument find an overtwisted disk when $\Sigma$ is a surface of 
genus greater than one and when the image of $B_x\subset \Sigma$ has a boundary component.
\ehw

So we can find an overtwisted disk unless there is some positive hyperbolic singularity on $\partial B_x$
(or $\Sigma=S^2$ and $e(\xi)([\Sigma])=0$). Thus we can cancel $x$ against a hyperbolic point. 
Continuing in this way we eventually show that $e(\xi)([\Sigma])\leq -\chi(\Sigma).$ (Note you should
be careful since when canceling $x$ new closed leaf may be born. If this happens add another pair
of negative singularities to break this closed leaf.)

\bhw
Finish the proof by showing that $-e(\xi)([\Sigma])\leq -\chi(\Sigma).$ Note you can do this by showing 
that $\Sigma$ may be perturbed so that $e_-=0.$ (Why is this sufficient?)
\ehw
\end{proof}

\bhw
Using the ideas in the proof of Theorem~\ref{genusbound} show: If there is an embedded disk $D$
in $(M,\xi)$ such that $D_\xi$ contains a closed leaf, then $\xi$ is overtwisted.  
The original definition of \dfn{tight} was the absence of embedded disks $D$ whose characteristic
foliation contains closed leaves. So it was not clear that a contact structure must be tight or
overtwisted. But this exercise shows that the original definition of tight
is equivalent to not being overtwisted.
\ehw

\section{Legendrian and Transverse Knots}\label{knots}

Just as studying surfaces in a contact 3-manifolds can illuminate the contact structure
so can studying curves. Two particularly interesting types of curves to study are
Legendrian curves and transverse curves. If $(M,\xi)$ is a contact manifold then a
curve $\gamma\co S^1\to M$ is called \dfn{Legendrian} (respectively \dfn{transverse}) if 
$\gamma(S^1)$ is always tangent (respectively transverse) to $\xi,$ that is for every
$x\in S^1,$ $d\gamma(T_xS^1)$ is contained in (respectively, is transverse to) $\xi_{\gamma(x)}.$
As is the custom in knot theory, we will frequently confuse $\gamma$ with its image. When we
try to classify Legendrian or transverse knots we will always be trying to classify them up
to isotopies through knots of the same type.

\begin{figure}[ht]
	{\epsfxsize=5truein\centerline{\epsfbox[0 737 394 842]{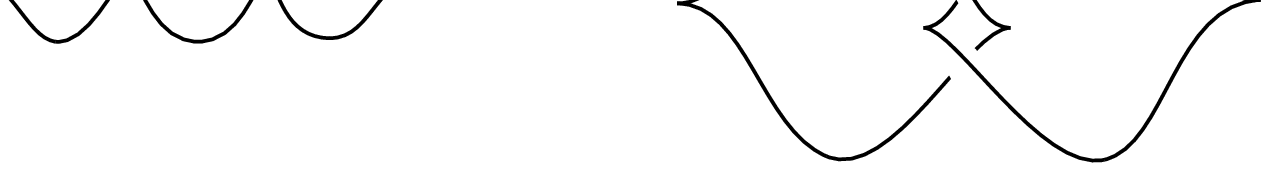}}}
	\caption{Examples of Legendrian knots.}
	\label{front}
\end{figure}

Let's begin by considering Legendrian and transverse knots in the standard contact structure 
on $\R^3.$ Recall, the contact structure is $\xi=\ker (dz+xdy).$ Now suppose $\gamma$ is
a Legendrian curve in $(\R^3, \xi).$ To picture $\gamma$ we will project it to the $yz$-plane.
This is called the \dfn{front projection} of $\gamma.$ The projection of $\gamma$ will ``look like''
Figure~\ref{front}. 
What we mean by ``look like'' is two things:
\begin{enumerate}
	\item at all the crossings the strand of $\gamma$ with the smaller 
	slope lies in front of the strand with the larger slope, and
	\item there are no vertical tangencies; instead there are cusps.
\end{enumerate}

\bhw
If $\gamma$ is Legendrian then show that 
	\begin{equation}\label{xcoord}
	x=-\frac{dz}{dy}.
	\end{equation}
	That is the $x$-coordinate of $\gamma$ is determined by the slope of its front
	projection. Thus the Legendrian knot $\gamma$ can be recovered from its front projection.
\ehw
\bhw
Convince your self that the restrictions above on the front projection are
	the only restrictions on a Legendrian knot and that any projection satisfying these
	restrictions is the projection of a Legendrian knot.

\begin{figure}[ht]
	{\epsfxsize=3.5truein\centerline{\epsfbox[0 389 492 842]{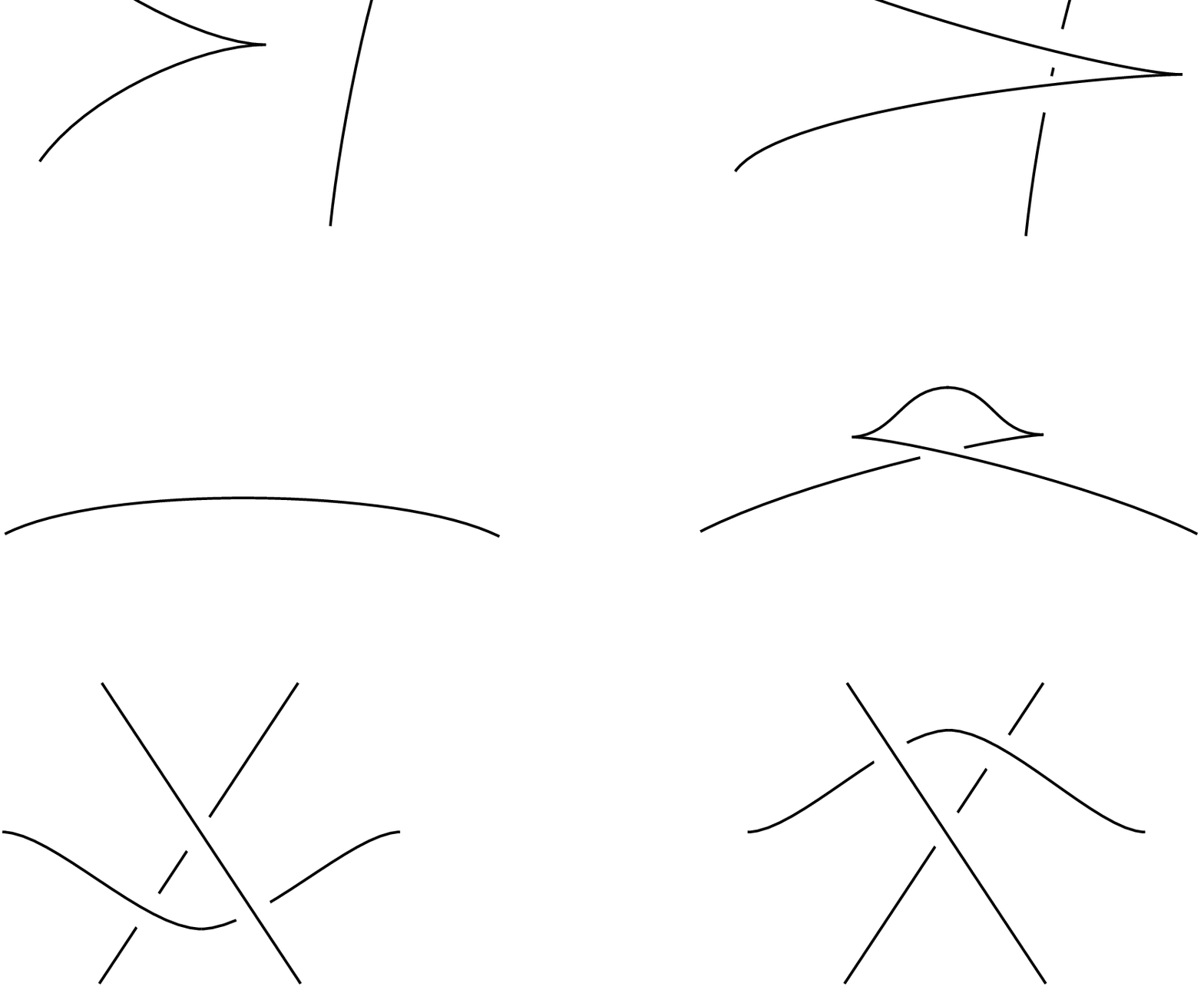}}}
	\caption{Legendrian Reidemister moves.}
	\label{moves}
\end{figure}

\ehw
From this exercise we see that the study of Legendrian knots in $\R^3$ reduces to the study of
their front projections. In particular if two Legendrian knots are isotopic (through Legendrian knots)
then you can get from the front projection of one to the front projection of the other by a
sequence of Legendrian Reidemister moves shown in Figure~\ref{moves} (and the moves obtained from these
by rotating the pictures $180^\circ$ around the $y$ or $z$-axes).
\begin{lem}\label{legapproxR}
	Any knot in $\R^3$ can be $C^0$ approximated by a Legendrian knot.
\end{lem}
\bhw
Prove this lemma.  \hfill\break
Hint: Consider the projection of the knot into the $yz$-plane. You will have
to consider how the projection fails to satisfy the two condition discussed above for a front projection
and how it fails to satisfy Equation~\eqref{xcoord}.
To fix this you can use ``zig-zags.'' 
\ehw

Even though we have been only discussing knots in $(\R^3,\xi)$ we can actually use this and
Darboux's Theorem to show
\begin{lem}\label{legapprox}
	Any curve in a contact manifold may be $C^0$ approximated by a Legendrian curve.
\end{lem}
\bhw
Prove this lemma.
\ehw

\bhw
Try to carry out a discussion, similar to the one above, for transverse knots. In other words
understand their ``front projections'' and prove the transverse versions of Lemmas~\ref{legapproxR} and 
\ref{legapprox}. If you get stuck take a look at \cite{EtnyreHonda??}.
\ehw

\subsection{The Classical Invariants of Legendrian and Transverse Knots}

The first step in trying to classify something is to find invariants that can help you distinguish
the objects under consideration ({\em e.g.} the Euler characteristic for surfaces). For Legendrian 
knots there are two easily defined invariants. Let $\gamma$ be a Legendrian knot and $\Sigma$ a surface
bounded by it. (If no such surface exists the situation is a bit more complicated, see \cite{EtnyreHonda01}.)
Take a vector field $v$ along $\gamma$ that is transverse to $\xi,$ then form $\gamma'$ by pushing
$\gamma$ in the direction of $v.$ Now the \dfn{Thurston-Bennequin invariant} of $\gamma,$
$\tb(\gamma),$ is the signed
intersection number of $\gamma'$ with $\Sigma$ ({\em i.e.} the linking number of $\gamma$ and $\gamma'$).
If we orient $\gamma$ then we can take a vector field $u$ along $\gamma$ that induces the chosen orientation
on $\gamma.$ Note that $u$ is in $\xi$ (since $\gamma$ is Legendrian). The \dfn{rotation number} of $\gamma,$
$r(\gamma),$
is Euler number $\xi|_\Sigma$ relative to $u.$ By this we mean $r(\gamma)$ is the obstruction to extending
$u$ to a non-zero vector field in $\xi|_\Sigma.$
\bhw
Choose any trivialization of $\xi$ over $\Sigma$. (Why can you always find such a trivialization?) 
Using this trivialization $u$ rotates some number of times as we traverse $\gamma$ positively ({\em i.e} in
the direction of the orientation). Prove that this number of rotations is the rotation number of $\gamma.$
\ehw

\begin{figure}[ht]
	{\epsfxsize=3.5truein\centerline{\epsfbox{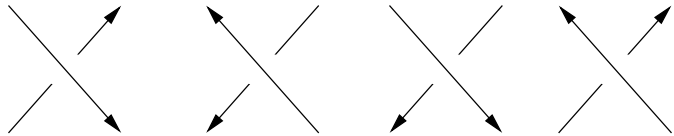}}}
	\caption{Right handed crossings (left) contribute $+1$ to the writhe while left
	handed crossings (right) contribute $-1.$}
	\label{signs}
\end{figure}

It is easy to compute these invariants in $\R^3$ using the front projection. Let $\gamma$ be an
oriented Legendrian knot in $\R^3.$ Recall the \dfn{writhe} of a knot diagram is the sum (over the
crossings in a diagram) of a $\pm 1$ at
each crossing, where the sign of the crossing is determined by its handedness. See Figure~\ref{signs}.
Denote  by $w(\gamma)$ the writhe of the front projection of $\gamma.$
Let $c(\gamma), c_u(\gamma)$ and $c_d(\gamma)$ be the number of cusps, upward oriented cusps and
downward oriented cusp (respectively) in the front projection. 
\begin{lem}
	With the notation above
	\begin{equation}
		\tb(\gamma)=w(\gamma)-\frac12 c(\gamma)
	\end{equation}
	and
	\begin{equation}
		r(\gamma)=\frac12(c_d(\gamma)-c_u(\gamma)).
	\end{equation}
\end{lem}
\bhw
Prove this Lemma.\hfill\break
Hint: A global trivialization of $\xi$ is given by $\{\frac{\partial}{\partial x},
x\frac{\partial}{\partial z}-\frac{\partial}{\partial y}\}.$ Use this trivialization to
compute the rotation number. Moreover,
the writhe of a diagram is the difference between the ``blackboard'' framing of a knot 
({\em i.e.} the obvious one coming from the diagram) and the framing coming from a Seifert surface.
Now use the vector $\frac{\partial}{\partial z}$ to compute the Thurston-Bennequin invariant.
\ehw
\bhw
Show $\tb$ and $r$ are invariants of Legendrian knots in $\R^3$ by using the Legendrian 
Reidemister moves. 
\ehw

Now suppose $\gamma$ is a transverse knot with Seifert surface $\Sigma.$ 
We choose a nonzero vector field
$v$ in $\xi|_\Sigma$ and form a copy of $\gamma'$ of $\gamma$ by pushing $\gamma$ in the direction
of $v.$ The \dfn{self-linking number} of $\gamma$ is the signed intersection number of $\gamma'$
with $\Sigma$ (once again it is just the linking number of $\gamma$ and $\gamma'$). 
The self-linking number of a knot in $\R^3$ may also be computed via its projection onto the
$yz$-plane. Specifically, one can show
\begin{equation}
	l(\gamma)=w(\gamma).
\end{equation}

\subsection{The Bennequin Inequality}

We may now state the fundamental Bennequin Inequality.
\begin{thm}
	If $\gamma$ is a transverse knot in a tight contact structure then
	\begin{equation}\label{bennequin}
		l(\gamma)\leq -\chi(\Sigma),
	\end{equation}
	where $\Sigma$ is any Seifert surface for $\gamma.$
\end{thm}
\bhw
Prove this theorem. \hfill\break
Hint: Since $\gamma$ is oriented it induces an orientation on its Seifert surface
$\Sigma.$ With these orientations the characteristic foliation
is oriented so that, thought of as a flow, it flows transversely out of $\partial \Sigma=\gamma.$
Thus $\partial \Sigma$ ``acts like a negative elliptic point.'' With this observation the proof
of this theorem is very similar to the proof of Theorem~\ref{genusbound}. It will be helpful to interpret
$l(\gamma)$ as a relative Euler class and then show (using notation from the proof of 
Theorem~\ref{genusbound}) that $l(\gamma)=-[(e_+-h_+)-(e_--h_-)].$
\ehw

This inequality provides a lower bound on the genus of a Seifert surface for $\gamma.$ 
In general, it is difficult to determine the smallest possible genus of a Seifert surface
for a given knot. (You should convince yourself that you can always find Seifert surfaces
of arbitrarily large genus for a given knot.) The Bennequin inequality can sometimes help in 
determining this smallest genus. 
\bhw
Look at the table of knots in \cite{Rolfson} and see which of them have
transverse realizations realizing the upper bound in Inequality~\eqref{bennequin}.
\hfill\break
Hint: It might be easier to consider Legendrian knots (see below). You will need to be able to construct
Seifert surfaces for the knots. The most common algorithm for this can be found in \cite{Rolfson}. 
\ehw

It is interesting to note that Bennequin proved Inequality~\eqref{bennequin} for any transverse
knot in the standard contact structure on $\R^3.$ But he did it {\em without knowing that the
contact structure was tight!} This, in fact, was the first hint that there was more than one
type of contact structure, but it still took several years for the notions of ``tight'' and ``overtwisted''
to be developed. So, in modern language, Bennequin proved the standard contact structure on
$\R^3$ was tight by proving Inequality~\eqref{bennequin}. Indeed, being able to prove this
inequality for a contact structure is equivalent to showing it is tight.
\bhw
Prove a contact structure is tight if and only if Inequality~\eqref{bennequin} is true. It might
be better to read below about the Legendrian version of Bennequin's Inequality and then 
think in terms of Legendrian knots.
\ehw
\noindent
He did this by examining relations between transverse knots and braid theory. See \cite{Birman} for
more on this relationship. Nowadays, our understanding of the inequality is somewhat different.
In stead of using it to prove a contact structure is tight, we usually prove the contact structure
is tight using other techniques and then use the inequality to study transverse knots in the
contact structure. The current preferred method to show contact structures are tight is to use
Theorem~\ref{fill2tight}

We now consider the Legendrian version of Bennequin's Inequality.
\begin{thm}\label{bennequin2}
	Let $\gamma$ be a Legendrian knot in a tight contact structure. Then
	\begin{equation}
		\tb(\gamma)+|r(\gamma)|\leq -\chi(\Sigma),
	\end{equation}
	where $\Sigma$ is a Seifert surface for $\gamma.$
\end{thm}
To prove this we just need to notice a simple relation between Legendrian and transverse
knots. Let $\gamma$ be a Legendrian knot and $A=[-\epsilon,\epsilon]\times S^1$ an embedded
annulus with $\gamma=\{0\}\times S^1$ and twisting so as never to be tangent to $\xi$ along $\gamma.$ 
Note that $\gamma$ is a closed leaf in $A_\xi$ and
for a generic choice of $A$ there will be no singularities and no other closed orbits. In fact,
we can assume (Why?)
that away from $\{0\}\times S^1$ the curves $\gamma_\pm=\{\pm x\}\times S^1,$ where 
$0<x<\epsilon,$ are transverse to $\xi.$ Furthermore, one can show that
\begin{equation}\label{pushoff}
	l(\gamma_\pm)=\tb(\gamma)\mp r(\gamma).
\end{equation}
\bhw
Understand the relation between $\gamma$ and $\gamma_\pm$ in the front projection. Use
this to prove Equation~\eqref{pushoff} for knots in the standard contact structure on $\R^3.$
You might want to try to prove the equation in general, or see \cite{EtnyreHonda??}.
\ehw

\bhw
Show how Theorem~\ref{bennequin2} follows from Equations~\eqref{bennequin} and \eqref{pushoff}.
\ehw

We have seen that the study of Legendrian and transverse knots can illuminate the nature of 
contact structures (such as the tight vs.\ overtwisted dichotomy), but their study is also quite
interesting in its own right. Legendrian and transverse 
unknots \cite{EliashbergFraser}, torus knots and figure eight knots \cite{EtnyreHonda??} have 
been classified and are essentially determined by their knot types and the invariants described
above. However, there are Legendrian knots that are topologically isotopic, have the same Thurston-Bennequin
invariants and rotation numbers but are not Legendrian isotopic. 
Such examples were first found in a tight contact structure on $S^2\times S^1\# S^2\times S^1$ (see
\cite{Fraser}). Here a geometric argument very specific to the situation was used to distinguish the
knots. Shortly after these examples were found an exciting new invariant was discovered
\cite{Chekanov, EGH, EtnyreNgSabloff} that allowed one to find many such ``non-simple'' Legendrian knots
in the standard tight contact structure on $S^3.$
The situation for transverse knots is not so well understood: it is unknown whether transverse
knots are determined by their topological knot type and their self-linking number.

\subsection{Transverse Knots and the Existence of Contact structures}

Dehn surgery is an important tool in understanding topological 3-manifolds. We wish to show that
Dehn surgery can be used in the world of contact 3-manifolds too. First let us recall the relevant 
definition. If $\gamma$ is a knot in a 3-manifold $M$ then it has a neighborhood, $N,$
diffeomorphic to $S^1\times D^2.$ Fix an embedded curve $\alpha$ on 
$\partial N\subset \partial \overline{M\setminus N}.$ Now choose any diffeomorphism $f$ of 
$T^2=\partial (S^1\times D^2)$ that sends the meridian, $\{p\}\times \partial D^2,$ to $\alpha$ and
define the \dfn{$\alpha$ Dehn surgery along $\gamma$} to be the manifold obtained from
$\overline{M\setminus N}$ by gluing in a solid torus via $f:$
\begin{equation}
	M(\gamma,\alpha)=(\overline{M\setminus N})\cup_f (S^1\times D^2).
\end{equation}
\bhw
Show that any choice of $f$ sending the meridian to $\alpha$ will produce the same 3-manifold
(up to diffeomorphism).
\ehw

We would now like to consider doing Dehn surgery on a transverse knot. To this end we observe that
another application of Moser's method yields
\begin{lem}\label{transnbhd}
	Let $\gamma_i$ be a transverse knot in $(M_i,\xi_i)$ for $i=1,2.$ Then any smooth 
	map from $\gamma_1$ to $\gamma_2$ may be extended to a contactomorphism from a neighborhood
	of $\gamma_1$ to a neighborhood of $\gamma_2.$
\end{lem}
Let's construct a standard model for the neighborhood of a transverse curve. For this
consider the contact structure $\xi=\ker(\cos r d\phi+r\sin r d\theta)$ on $S^1\times \R^2,$ where $\phi$ is
the coordinate on $S^1$ and $(r,\theta)$ are polar coordinates on $\R^2.$ (Note this contact structure
is just the one in Example~\ref{otexample} 
with the $z$-axis wrapped around the $S^1.$ Said another way, $\R^3$ is
the universal cover of $S^1\times\R^2$ and the contact structure in Example~\ref{otexample} is just the
pull back of this one under the covering map.) Note that $T_a=\{(\phi, r, \theta) | r=a\}$ is a
torus, and $(T_a)_\xi$ is a non-singular foliation by lines of slope $-a\tan a.$ 
Lemma~\ref{transnbhd}
implies that any transverse knot $\gamma$ has a neighborhood $N$ contactomorphic to 
$S_a=\{(\phi, r, \theta) | r\leq a\}$ for some $a.$ 

Now if $\gamma$ is some transverse knot in $S^3,$ with the standard contact structure, 
then it has a neighborhood $N$ contactomorphic to $S_a$ for some $a.$ If we remove $N$ from
$S^3$ and then glue in a solid torus $S^1\times D^2$ via a map $f,$ the resulting manifold
$M$ has a contact structure defined on all but the $S^1\times D^2$ part. Note that on $\partial(S^1\times
D^2)\subset M$ we have a characteristic foliation. This foliation is a linear foliation with
some slope $s$ (when measured with respect to the $S^1\times D^2$ product structure).
\bhw
Determine what $s$ is in terms of $a$ and the slope of the curve $\alpha.$ Is $s$ uniquely determined?
If not what are the possible $s$'s.
\ehw
Now we can find a model contact structure on $S^1\times D^2$ whose characteristic foliation is
also linear with slope $s.$ 
\bhw
Check that this model contact structure on $S^1\times D^2$ and the contact structure on 
$\overline{S^3\setminus N}$ induced from $S^3$ define a contact structure on $M.$
\ehw
We can clearly perform this construction on a link in $S^3.$ Thus since any 3-manifold can be
obtained from $S^3$ by Dehn surgery on a link we have proved:
\begin{thm}[Martinet \cite{Mar71}]
	All closed compact 3-manifolds support a contact structure.
\end{thm}
Note that there are many choices for $a$ so that $S_a$ has the appropriate slope to be
used in the above construction. However, it is clear that if we choose any $a$ except the
smallest possible $a$ then we automatically get an overtwisted structure. (Find the overtwisted
disk!) Even choosing the smallest possible $a$ we will frequently get an overtwisted structure
on the surgered manifold. If you are sufficiently careful with this construction you can show
\begin{thm}[Lutz \cite{l}]
	In every homotopy class of oriented plane fields on a closed compact 3-manifold
	there is an overtwisted contact structure.
\end{thm}
\bhw
Try to prove this theorem on $S^3.$\hfill\break
Hint: There are $\Z$ homotopy classes of oriented plane fields. (To see
this trivialize the tangent bundle and choose a metric. Now given a plane field you can use
the unit vector orthogonal to the plane field to get a map to $S^2,$ well defined up to homotopy.
Thus homotopy classes of plane fields are in one-to-one correspondence with homotopy classes
of maps $S^3\to S^2.$ That is $\pi_3(S^2)=\Z.$) The standard contact structure on $S^3$ is orthogonal
to the Hopf fibration of $S^3.$  So if $\gamma$ is a fiber in the Hopf fibration then it has a neighborhood
contactomorphic to $S_a$ for some $a.$ Now replace $S_a$ with $S_{a+b},$ where $b$ is chosen so that
$S_a$ and $S_{a+b}$ have characteristic foliations with the same slope. Note that when we do this
we are still on $S^3$ but the contact structure is (possibly) different. 
\ehw

Unfortunately it is much harder to construct tight contact structures. 
\bhw
Try to understand how the constructions above relate to ``Legendrian surgery'' in which tight contact
structures are produced. See \cite{Gompf98} for a discussion of Legendrian surgery and
\cite{Gay} for part of its relation to the surgery described above.
\ehw

\section{Introduction to Convex Surfaces}

In the previous sections we have been discussing a classical approach to contact geometry.
By classical, I mean concentrating on specific characteristic foliations. In \cite{Giroux91}, 
Giroux initiated the use of convex surfaces in contact geometry. Using the theory of convex surfaces
one can ignore specific characteristic foliations when studying surfaces in a contact structure and 
concentrate on a few curves on the surface (the so called ``dividing curves''). In this section we will
indicate how to use convex surfaces in the study of contact geometry. For applications of this to
the classification of contact structures see \cite{Honda1} and \cite{Giroux00}, to the classification
of Legendrian knots see \cite{EtnyreHonda??} and to the nature of tightness see \cite{EtnyreHonda01}.

Given a contact manifold $(M,\xi)$ a vector field is called \dfn{contact} if its flow preserves the 
contact structure. A surface $\Sigma$ is called \dfn{convex} if there is a contact vector field 
transverse to it.
\bhw\label{vertinvt}
Show a surface $\Sigma$ is convex if and only if there is a neighborhood $N=\Sigma \times I$ such
that $\xi|_N$ is invariant in the $I$ direction. (Note this exercise implies that convex is a not
such a great term for such a surface but we are stuck with it.)
\ehw
The first question one should ask is: Are there any convex surfaces? In \cite{Giroux91} it was shown 
that any closed surface is $C^\infty$-close to a convex surface. Moreover, in \cite{Kanda95} it was shown 
that this is also true for a surface with boundary so long as the surface has Legendrian boundary
and the twisting of the contact planes relative to the surface is not positive.

Now let $\Gamma$ be the set of points on a convex surface $\Sigma$ where the contact vector field
is tangent to $\xi.$  In \cite{Giroux91} it was shown that
(generically)
$\Gamma$ is a multi-curve, that is collection of curves, on $\Sigma.$ The multi-curve 
$\Gamma$ satisfies:
\begin{enumerate}
\item $\Sigma\setminus\Gamma=\Sigma_+\coprod \Sigma_-,$
\item $\Sigma_\xi$ is transverse to $\Gamma,$ and
\item there is a vector field $w$ and volume form $\omega$ on $\Sigma$ such that
\begin{enumerate}
\item $w$ is directs $\Sigma_\xi$ ({\em i.e.} $w$ is tangent to $\Sigma_\xi$ where it is nonsingular and
	is zero where it is singular),
\item the flow of $w$ expands $\omega$ on $\Sigma_+$ and contracts $\omega$ on $\Sigma_-$,
\item $w$ points transversely out of $\Sigma_+.$
\end{enumerate}
\end{enumerate}
\bhw
Verify these properties for $\Gamma.$\hfill\break
Hint: Use Exercise~\ref{vertinvt} to show that in a neighborhood of $\Sigma,$ $\xi$ is the kernel
of $\beta+u dt$ where $\beta$ is a 1--form on $\Sigma$ and $u$ is a function on $\Sigma.$ Now try
to understand $\Gamma$ and $\Sigma_\xi$ in terms of this 1--form.
\ehw
If $\mathcal{F}$ is any singular foliation of $\Sigma$ then a multi-curve $\Gamma$ on $\Sigma$ is
said to \dfn{divide} $\mathcal{F}$ if they satisfy the above conditions where $\Sigma_\xi$ is replaced
by $\mathcal{F}.$ Moreover, the curves $\Gamma$ are called 
the \dfn{dividing curves} for  $\mathcal{F}.$ We now have
the first major theorem about convex surfaces.
\begin{thm}[Giroux \cite{Giroux91}]\label{curvesdet}
	Suppose $\mathcal{F}$and $\Sigma_\xi$ are both divided by the same multi-curve $\Gamma.$ Then
	inside any neighborhood $N$ of $\Sigma$ there is an isotopy $\Phi_t\co \Sigma\to N, t\in[0,1]$  of 
	$\Sigma$  such that
	\begin{enumerate}
		\item $\Phi_0= \mbox{inclusion of $\Sigma$ into $N$},$
		\item $\Phi_t(\Sigma)$ is a convex surface for all $t,$
		\item $\Phi_t$ does not move $\Gamma,$
		\item $(\Phi_1(\Sigma))_\xi=\Phi_1(\mathcal{F}).$
	\end{enumerate}
\end{thm}
\noindent
This theorem basically says that given a convex surface we can assume the characteristic foliation
is anything we wish it to be as long as it is divided by the appropriate curves. Or said another way,
it is really the dividing curves that carry the essential information about the contact structure
in a neighborhood of a convex surface and not the 
specific characteristic foliation. One needs to be very careful with this
heuristic statement but it is a useful way to think about convex surfaces. 
\bhw\label{otdiskfromdivides}
Show that if $\Sigma\not =S^2$ is any convex surface in $(M,\xi)$ and it has a closed contractible dividing
curve then $\xi$ is overtwisted. \hfill\break
Hint: Try to write down some foliation respecting the dividing curves in which it is easy to
see an overtwisted
disk. Why is it important that $\Sigma\not=S^2?$
\ehw

\begin{figure}[ht]
	{\epsfysize=1.5truein\centerline{\epsfbox[0 701 226 841]{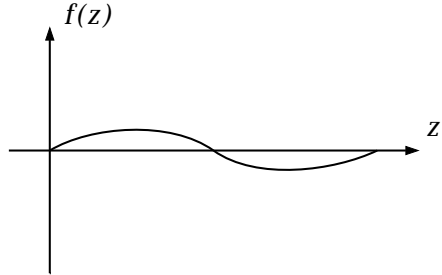}}}
	\caption{The graph of $f.$}
	\label{graph}
\end{figure}

\bex
Consider $\R^3$ with the contact structure $\xi=\ker(dz+xdy).$ If we quotient $\R^3$ by $z\mapsto z+1$
and $y\mapsto y+1$ we will get $M=\R\times T^2$ and since the contact structure is preserved by this
action $\xi$ will induce a contact structure on $M.$ The characteristic foliation on $\{0\}\times T^2$
is by horizontal lines ({\em i.e.} by the lines $z=\hbox{constant}$). You can check that this is
{\em not} a convex surface, but it is easy to perturb into a convex surface. Let $f\co [0,1]\to \R$ be the
function whose graph is given in Figure~\ref{graph},
then set $\Sigma=\{(f(z), y,z)\}.$
\begin{figure}[ht]
	{\epsfysize=1.5truein\centerline{\epsfbox[0 698 349 841]{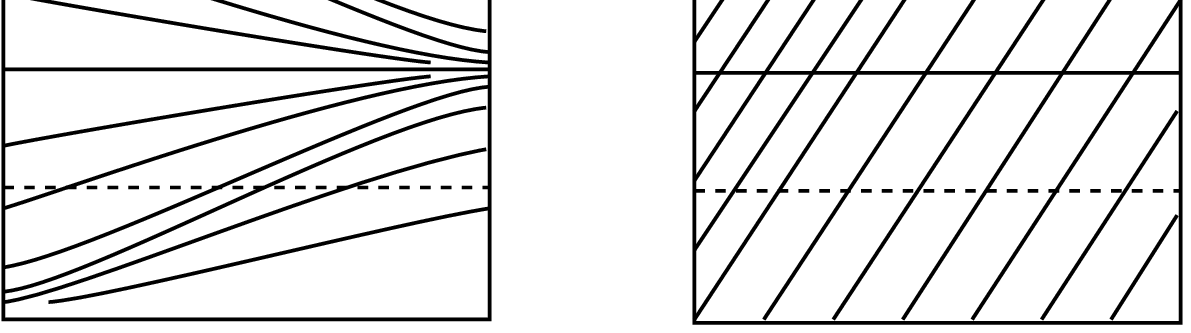}}}
	\caption{The characteristic foliation on $\Sigma$ (left), with
	dividing curves (dotted lines). Another foliation on $\Sigma$ with the same dividing set (right).}
	\label{convexex}
\end{figure}
Clearly $\Sigma$ is a small perturbation of $T^2\times\{0\},$ moreover, the characteristic foliation
on $\Sigma$ is as shown on the left hand side of Figure~\ref{convexex}. It is easy to check
that the dotted lines in the figure give a set of dividing curves for $\Sigma_\xi,$ and thus $\Sigma$
is convex. Now using Theorem~\ref{curvesdet} we can perturb $\Sigma$ so as to realize any 
characteristic foliation that respects these dividing curves. In particular, we can arrange for the foliation
to look like the one the right hand side of Figure~\ref{convexex}. This foliation has two lines of
singularities along $\{z=0\}\cup\{z=\frac12\}$ and all of the nonsingular leaves have slope $s\not=0.$
The nonsingular leaves are called \dfn{ruling curves} and the singular curves are called \dfn{Legendrian
divides}. Note the Legendrian divides must be parallel to the dividing curves, but we may choose the
ruling curves to have any slope except 0, the slope of the Legendrian divides. Any torus with 
a foliation like this will be said to be in \dfn{standard form}. 
\eex

The power of convex surfaces is contained largely in Theorem~\ref{curvesdet} in conjunction with 
\begin{lem}[\cite{Kanda95, Honda1}]\label{lem:trans}
	Suppose that $\Sigma$ and $\Sigma'$ are convex surfaces, with dividing curves $\Gamma$ and $\Gamma',$
	and $\partial \Sigma'\subset \Sigma$ is Legendrian. 
	Let $S=\Gamma\cap \partial\Sigma'$ and $S'=\Gamma'\cap \partial\Sigma'.$ 
	Then between each two adjacent points in $S$ there is one point in $S'$ and vice verse.
	See Figure~\ref{transfer}.
\begin{figure}[ht]
	{\epsfysize=2.2truein\centerline{\epsfbox[0 519 300 840]{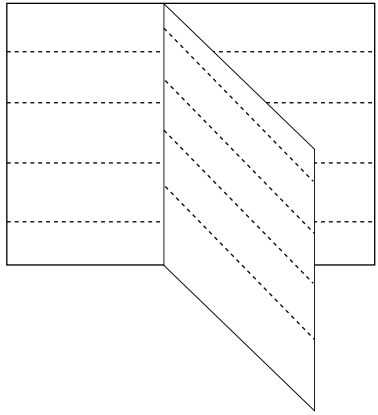}}}
	\caption{Transferring information about dividing curves from one surface to another. The top
	and bottom of the picture are identified.}
	\label{transfer}
\end{figure}
	(Note the sets $S$ and $S'$ are cyclically ordered since they sit on $\partial\Sigma'$)
\end{lem}
In this lemma clearly $\Sigma'$ is not a closed surface. All of our previous discussion goes through
for surface with boundary as long as the boundary is Legendrian and the twisting of the contact
planes relative to the surface in not positive. See \cite{Kanda95}.

We can now give a simple proof of the following result which is essentially due to Makar-Limanov \cite{ml},
but for the form presented here see Kanda \cite{Kanda95}. Though this theorem
seems easy, it has vast generalizations which we indicate below.
\begin{thm}\label{ctor} 
Suppose $M=D^2\times S^1$ and $\mathcal{F}$ is a singular foliation on $\partial M$ that is
divided by two parallel curves with slope $\frac1n$ (here slope $\frac1n$ means that the curves
are homotopic to $n[\partial D^2\times\{p\}]+[\{q\}\times S^1]$ where $p\in S^1$ and $q\in \partial D^2$).
Then there is a unique tight contact structure on $M$ whose characteristic foliation on $\partial M$ 
is $\mathcal{F}.$
\end{thm}
\begin{proof}
Suppose we have two tight contact structures $\xi_0$ and $\xi_1$ on $M$ inducing $\mathcal{F}$ as the
characteristic foliation on $\partial M.$  We will find a contactomorphism from $\xi_0$ to $\xi_1$
(in fact this contactomorphism will be isotopic to the identity). Let $f\co  M \to  M$ 
be the identity map. By Theorem~\ref{foldet} we can isotop $f$ rel.~$\partial M$ 
to be a contactomorphism in
a neighborhood $N$ of $\partial M.$  Now let $T$ be a convex torus in $N$ isotopic to $\partial M.$ Moreover
we can assume that the characteristic foliation on $T$ is in standard form. We know the slope of the
Legendrian divides is $\frac1n$ and we choose the slope of the ruling curves to be 0. Let $D$ be
a meridianal disk whose boundary is a ruling curve. We can perturb $D$ so that it is convex and
using Lemma~\ref{lem:trans} we know that the dividing curves for $D$ intersect the boundary of
$D$ in two points. Moreover, since there are no closed dividing curves on $D$ (since the contact
structure is tight, see Exercise~\ref{otdiskfromdivides}) 
we know that $\Gamma_D$ consists of one arc. We may isotop $f(D)$ (rel.~boundary)
to $D'$ so that all of this is true for $D'$ with respect to $\xi_1.$ 
Now using Theorem~\ref{curvesdet} we can arrange
that the characteristic foliations on $D$ and $D'$ agree; and further, we can isotop $f$ (rel.~$N$) so
that $f$ takes $D$ to $D'$ and preserves the characteristic foliation on $D.$ Thus another application
of Theorem~\ref{foldet} says we can isotop $f$ so as to be a contactomorphism on $N'=N\cup U,$ where
$U$ is a neighborhood of $D.$  Note that $B=\overline{M\setminus N'}$ is a 3--ball, so Theorem~\ref{3ball}
tells us that we can isotop $f$ on $B$ so that it is a contactomorphism on $B$ too. Thus $f$ is
a contactomorphism on all of $M$ and we are done with the proof. 
\end{proof}

\bhw\label{torus}
Suppose that $\mathcal{F}$ is a convex foliation on $\partial(S^1\times D^2)$ 
with $2n$ dividing curves of slope $\frac{p}{q}.$ Find an upper bound on the number
of tight contact structures on $S^1\times D^2$ which induce this foliation. For $\frac{p}{q}=\frac2q$ 
or $\frac3q$ (and $n=1$) classify the corresponding tight contact structures. If you are feeling
bold you might want to try and prove the upper bound you found is not in general sharp and then actually
find the sharp upper bound. This second part is not particularly easy; if you would like to see 
the answer consult \cite{Honda1}.
\ehw

\bhw
Try to generalize Exercise~\ref{torus} to a genus $g$-handle body. Which configurations
of dividing curves correspond to a unique contact structure? Can their ever be infinitely
many tight structures with a fixed foliation?
\ehw

\bhw
Prove the well known folk theorem of Eliashberg: There is a unique positive tight contact structure on
$S^1\times S^2.$ \hfill\break 
Hint: The argument is similar to the one in the proof of Theorem~\ref{ctor}. First
see that you can normalize the contact structure in a neighborhood of $\{p\}\times S^2.$ The
complement of this neighborhood is $I\times S^2,$ where $I$ is an interval. Now normalize the contact
structure in the neighborhood of an annulus $A= I\times S^1,$ for an appropriately chosen $S^1\subset S^2.$
Warning: Be careful here, you need to find a way to deal with the fact that 
the dividing curves on $A$ can spin around the $S^1$ factor many times.
The complement of the normalized regions is a 3--ball which has a unique tight contact structure (with
given boundary data).
\ehw

These exercises only give a hint at the power of convex surfaces. For further developments see
\cite{EtnyreHonda01, EtnyreHonda??, Giroux00, Honda1, Honda2}.


\end{document}